\documentclass[12pt]{amsart}
\usepackage[title,toc,titletoc,header]{appendix}
\usepackage[english]{babel}
\usepackage[utf8x]{inputenc}
\usepackage[T1]{fontenc}
\usepackage{amsmath}
\usepackage{amstext}
\usepackage{amsfonts}
\usepackage{amssymb}
\usepackage{amsthm}
\usepackage{amsrefs}
\usepackage{mathtools}
\usepackage{dsfont}
\usepackage{color}
\usepackage{graphicx}
\usepackage[colorinlistoftodos]{todonotes}
\usepackage[colorlinks=true, allcolors=blue]{hyperref}
\usepackage{float}
\allowdisplaybreaks
\usepackage{enumitem}
\usepackage{lipsum}
\usepackage{graphicx}
\newtheorem{theorem}{Theorem}[section]
\newtheorem{lemma}[theorem]{Lemma}
\newtheorem{prop}[theorem]{Proposition}
\newtheorem{cor}[theorem]{Corollary}

\newtheorem*{cor*}{Corollary}
\newtheorem*{thm*}{Theorem}
\newtheorem*{lem*}{Lemma}

\newtheorem*{prop*}{Proposition}
\theoremstyle{definition}

\newtheorem{example}[theorem]{Example}

\newtheorem*{defn*}{Definition}

\theoremstyle{remark}

\newcommand{\cA}{\mathcal{A}}

\newcommand{\bE}{{\mathbb{E}}}

\newcommand{\thissubtheorem}{}
\newcommand{\subtheorem}[1][]{%
  \if\relax\detokenize{#1}\relax
    \def\thissubtheorem{}%
  \else
    \def\thissubtheorem{ (#1)}%
  \fi
  \item
}

\title{on Intermediate subalgebras of $C^*$-simple Group Actions}

\author[]{Tattwamasi Amrutam \\with an Appendix by Tattwmasi Amrutam and Yongle Jiang}

\address{Department of Mathematics\\ University of Houston\\USA}
\email{tamrutam@math.uh.edu}
\address{Institute of Mathematics Polish Academy of Sciences\\ Warsaw, Poland\\}
\email{yjiang@impan.pl}

\date{}

\begin{document}

\begin{abstract}
We show that for a large class of actions $\Gamma \curvearrowright \mathcal{A}$ of $C^*$-simple groups $\Gamma$ on unital $C^*$-algebras $\mathcal{A}$, including any non-faithful action of a hyperbolic group with trivial amenable radical, every intermediate $C^*$-subalgebra $\mathcal{B}$, $C_{\lambda}^*(\Gamma)\subseteq \mathcal{B} \subseteq \mathcal{A}\rtimes_{r}\Gamma$, is of the form $\mathcal{A}_1\rtimes_{r}\Gamma$,  where $\mathcal{A}_1$ is a unital $\Gamma$-$C^*$-subalgebra of $\mathcal{A}$. 
\end{abstract}

\maketitle

\section{Introduction and Statement of the main result}

This note is concerned with the structure of intermediate subalgebras i.e., $C^*$-subalgebras $\mathcal{B}$, $C_{\lambda}^*(\Gamma) \subseteq \mathcal{B} \subseteq \mathcal{A}\rtimes_r\Gamma$, where $\Gamma \curvearrowright \mathcal{A},$ is an action of a $C^*$-simple group $\Gamma$ on a unital $C^*$-algebra $ \mathcal{A}$. For a large class of such actions, we give complete description of intermediate subalgebras. Namely, we show that if the kernel of the action contains a subgroup with respect to which, $\Gamma$ has the Powers' averaging property, all the intermediate subalgebras are of the form of  crossed product.
\begin{defn*} A non-trivial $\Lambda \le \Gamma$ is called a \emph{plump subgroup} if the following holds: for every $\epsilon >0$ and finite subset $F \subset \Gamma\setminus\{e\}$, there exist $s_1,s_2,\ldots,s_m \in \Lambda$ such that 
$$\left\|\frac{1}{m}\sum_{j=1}^m\lambda_{s_j}\lambda_t\lambda_{s_j^{-1}}\right\| < \epsilon, ~~\forall t \in F.$$
\end{defn*}

 \begin{theorem}
\label{Main Theorem}
Let $\Gamma$ be a discrete group with the approximation property (AP), let $\mathcal{A}$ be a unital $\Gamma$-$C^*$-algebra. Suppose that the kernel of the action $\Gamma \curvearrowright \mathcal{A}$ contains a plump subgroup of $\Gamma$. Then, every intermediate $C^*$-subalgebra $\mathcal{B}$, $C_{\lambda}^*(\Gamma) \subseteq \mathcal{B} \subseteq \mathcal{A}\rtimes_{r}\Gamma$, is of the form $\mathcal{A}_1\rtimes_{r}\Gamma$, for some unital $\Gamma$-$C^*$-subalgebra $\mathcal{A}_1$ of $\mathcal{A}$. 
\end{theorem}

We give examples of several classes of actions $\Gamma \curvearrowright \mathcal{A}$ which fall into the premises of Theorem \ref{Main Theorem}, in the last section. In particular, we show the following.
 \begin{theorem}
\label{example}
Let $\Gamma$ be a hyperbolic group with trivial amenable radical. For any non-faithful action $\Gamma \curvearrowright \mathcal{A}$, every intermediate $C^*$-subalgebra $\mathcal{B}$, $C_{\lambda}^*(\Gamma) \subseteq \mathcal{B} \subseteq A\rtimes_{r}\Gamma$, is of the form $\mathcal{A}_1\rtimes_{r}\Gamma$, where $\mathcal{A}_1$ is a $\Gamma$-$C^*$-subalgebra of $\mathcal{A}$.
\end{theorem}In a recent paper \cite{Suz18}, Suzuki proved a similar result giving a complete description of intermediate subalgebras, for both $C^*$-algebras and von Neumann algebras. 

The paper also contains a section which includes our joint result with Yongle Jiang proving a similar result for Odometer actions. 
\section*{Acknowledgement}
 The author is grateful to Mehrdad Kalantar, who is his advisor, for his helpful comments and suggestions on the organization of this paper. We thank Yongle Jiang and Ilya Kapovich for their helpful comments, and Nikolaos Panagopoulos for pointing out typo in an earlier version of the paper.
 
\section{Proof of Main Theorem}
For a faithful $*$-representation $\pi: \cA \to \mathbb{B}(\mathcal{H})$ of $\cA$ into the space of bounded operators on the Hilbert space $\mathcal{H}$, the crossed product $C^*$-algebra $\mathcal{A}\rtimes_r\Gamma$ is generated (as a $C^*$-algebra inside $\mathbb{B}(\ell^2(\Gamma, \mathcal{H}))$, the space of square summable $\mathcal{H}$-valued functions on $\Gamma$) by a unitary representation $\lambda$ of $\Gamma$ and a faithful $*$-representation of $\mathcal{A}$ in $\mathbb{B}(\ell^2(\Gamma, \mathcal{H}))$. This representation translates the action $\Gamma \curvearrowright \mathcal{A}$ into an inner action by the unitaries $\lambda_s, s \in \Gamma$. 

The reduced crossed product $\mathcal{A}\rtimes_{r}\Gamma$ comes equipped with the faithful canonical conditional expectation, $\bE: \mathcal{A}\rtimes_{r}\Gamma \to\cA$, which is defined by $\bE(a) = a$ and $\bE(a\lambda_s) = 0$ for $a\in \cA$ and $s \in \Gamma\setminus\{e\}$. We refer the reader to \cite{BroOza08} for more details on this.

\hfill{}

\noindent
{\it Proof of Theorem \ref{Main Theorem}.}
\label{C-case}
Let $\mathcal{A}$ be a unital $\Gamma$-$C^*$-algebra and let $\mathcal{B}$ an intermediate $C^*$-subalgebra of the form $C_{\lambda}^*(\Gamma) \subseteq \mathcal{B} \subseteq \mathcal{A}\rtimes_r\Gamma$. Suppose that $\Lambda$ is a plump subgroup of $\Gamma$ such that $\Lambda$ is contained in the kernel of the action $\Gamma \curvearrowright \mathcal{A}$. Fix $b \in \mathcal{B}$. Let $\epsilon>0$. Then, there are $t_1,t_2,\ldots,t_n \in \Gamma\setminus\{e\}$ such that  \[\left\|b-\left(\sum_{i=1}^na_{t_i}\lambda_{t_i}+\mathbb{E}(b)\right)\right\|_{\mathbb{B}(l^2(\Gamma,\mathcal{H}))} < \epsilon.\] Let $M=\max_{1 \le i \le n}\|a_{t_i}\|_{\mathcal{A}}$. Since $\Lambda$ is a plump subgroup of $\Gamma$, there exist $s_1,s_2,\ldots,s_m \in \Lambda$ such that \[\left\|\frac{1}{m}\sum_{j=1}^m\lambda_{s_j}\lambda_t\lambda_{s_j^{-1}}\right\|_{\mathbb{B}(l^2(\Gamma))} < \frac{\epsilon}{nM}, ~~\forall i=1,2,\ldots,n.\]
By \cite[Lemma 2.1]{A-K}, it follows that 
  \[\begin{split}&\left\|\frac{1}{m}\sum_{j=1}^m\lambda_{s_j}\left(\sum_{i=1}^na_{t_i}\lambda_{t_i}\right)\lambda_{s_j^{-1}}\right\|_{\mathbb{B}(l^2(\Gamma,\mathcal{H}))} \\&\le \sum_{i=1}^n\|a_{t_i}\|_{\mathcal{A}}\left\|\frac{1}{m}\sum_{j=1}^m\lambda_{s_j}\lambda_t\lambda_{s_j^{-1}}\right\|_{\mathbb{B}(l^2(\Gamma))}  < \epsilon.\end{split}\]
 Now,
          \[ \begin{split} &\left\|\frac{1}{m}\sum_{j=1}^m\lambda_{s_j}(b-\bE(b))\lambda_{s_j^{-1}}\right\|_{\mathbb{B}(l^2(\Gamma,\mathcal{H}))}
     \\&\le\left\|\frac{1}{m}\sum_{j=1}^m\lambda_{s_j}\left(b-\left(\sum_{i=1}^na_{t_i}\lambda_{t_i}+\mathbb{E}(b)\right)\right)\lambda_{s_j^{-1}}\right\|_{\mathbb{B}(l^2(\Gamma,\mathcal{H}))} 
     \\&+\left\|\frac{1}{m}\sum_{j=1}^m\lambda_{s_j}\left(\sum_{i=1}^na_{t_i}\lambda_{t_i}\right)\lambda_{s_j^{-1}}\right\|_{\mathbb{B}(l^2(\Gamma,\mathcal{H}))}< 2\epsilon.
        \end{split}\]
Since $\Lambda$ acts trivially on $\mathcal{A}$, we get that 
\[\begin{split}&\left\|\frac{1}{m}\sum_{j=1}^m\lambda_{s_j}b\lambda_{s_j^{-1}}-\bE(b)\right\|_{\mathbb{B}(l^2(\Gamma,\mathcal{H}))}\\&=\left\|\frac{1}{m}\sum_{j=1}^m\lambda_{s_j}\left(b-\bE(b)\right)\lambda_{s_j^{-1}}\right\|_{\mathbb{B}(l^2(\Gamma,\mathcal{H}))}< 2\epsilon\end{split}\]
Since $\epsilon >0$ is arbitrary, this shows that $\bE(\mathcal{B})\subset \mathcal{B}$. By \cite[Proposition $3.4$]{Suz17}, it follows that $\mathcal{B}=\mathbb{E}(\mathcal{B})\rtimes_r\Gamma$. \qed

\hfill{}

A similar result holds in the von Neumann algebraic setting. We thank Yongle Jiang for pointing out to us that the same proof as above works in this setting as well. 

\hfill{}
\hfill{}

\noindent
\begin{theorem}
\label{VN}
Let $\Gamma$ be a discrete group, $\mathcal{M}$ be a $\Gamma$-von Neumann-algebra with a separable predual. Suppose that $\Lambda$ is a plump subgroup of $\Gamma$ such that $\Lambda$ is contained in the kernel of the action $\Gamma \curvearrowright \mathcal{M}$. Then every intermediate von Neumann-subalgebra $\mathcal{N}$, $L(\Gamma) \subseteq \mathcal{N} \subseteq \mathcal{M}\rtimes\Gamma$ is of the form $\mathcal{M}_1\rtimes\Gamma$, where $\mathcal{M}_1$ is a $\Gamma$-von Neumann subalgebra of $\mathcal{M}$. 
\label{VN-Case}
\begin{proof}
Let $\varphi$ be a faithful normal state on $\mathcal{M}$ and let $\tilde{\varphi}(a)=\varphi(\mathbb{E}(a)), a \in \mathcal{M}\rtimes\Gamma$, where $\mathbb{E}$ is the canonical conditional expectation from $\mathcal{M}\rtimes\Gamma$ onto $\mathcal{M}$. Then $\tilde{\varphi}$ is a faithful normal state on $\mathcal{M}\rtimes\Gamma$. Consider the $\|.\|_2$-norm on $\mathcal{M}\rtimes\Gamma$ associated to $\tilde{\varphi}$, defined by \[\|a\|_2:=\sqrt{\tilde{\varphi}(a^*a)} \text{ ~ for } a \in \mathcal{M}\rtimes\Gamma.\]
Let $b \in \mathcal{N}$ and let $\epsilon>0$ be given. Then there are $t_1,t_2,\ldots,t_n \in \Gamma\setminus{e}$ such that 
\[\left\|b-\left(\sum_{i=1}^na_{t_i}\lambda_{t_i}+a_e\right)\right\|_2 < \frac{\epsilon}{2}.\]
Since $\tilde{\varphi}$ is $\mathbb{E}$-invariant, $\mathbb{E}$ is continuous with respect to the $\|.\|_2$-norm and hence by the triangle inequality, we see that
\[\left\|b-\left(\sum_{i=1}^na_{t_i}\lambda_{t_i}+\mathbb{E}(b)\right)\right\|_2 <\epsilon.\]
Since the $\|.\|_2$ is dominated by the operator norm, by proceeding exactly as in the proof of Theorem \ref{Main Theorem}, we see that $\mathbb{E}(b) \in \mathcal{N}$. By \cite[Corollary 3.4]{Suz18}, the proof is complete. 
\end{proof}
\end{theorem}

\section{Examples}
In this section, we give examples of classes of actions of $\Gamma$ on unital $C^*$-algebra $\mathcal{A}$, which satisfy the conditions of Theorem \ref{Main Theorem}. 
\subsection{Plump Subgroups}
First, we give several conditions under which a subgroup $\Lambda$ is plump in $\Gamma$.  

\begin{theorem}
\label{lem}
Let $\Lambda$ be a subgroup of $\Gamma$. Suppose that there is a free action of $\Gamma$ on a compact Hausdorff space $X$ such that the action of $\Gamma$ restricted to $\Lambda$ is strongly proximal. Then $\Lambda$ is plump in $\Gamma$. 
 \begin{proof} First, proceeding exactly as in the proof of \cite[Theorem 4.5]{Haagerup}, one sees that $\tau_0 \in \overline{\{s\varphi: s \in \Lambda\}}^{\text{weak*}}$ for each state $\varphi$ on $C_{\lambda}^*(\Gamma)$. Again, by the same theorem, for any finite collection  $t_1,t_2, \ldots, t_n \in \Gamma-\{e\}$ and $\epsilon > 0$, there exist $s_1,s_2,\ldots,s_m \in \Lambda$ such that \[
\left\|\frac{1}{m}\sum_{j=1}^m\lambda_{s_j}\lambda_{t_i}\lambda_{s_j^{-1}}\right\| < \epsilon, \text{ for each $i=1,2,\ldots,n$ }.
\]
\end{proof}
\end{theorem}

 \begin{cor}
 \label{3.2}
 Suppose that $\Lambda$ is a normal $C^*$-simple subgroup of $\Gamma$ with trivial centralizer inside $\Gamma$.  Then $\Lambda$ is plump in $\Gamma$. 
 \begin{proof}
By {\cite[Lemma 5.3]{BKKO}}, the action $\Lambda \curvearrowright \partial_F\Lambda$ extends to a free action of $\Gamma$ on $ \partial_F\Lambda$. By Theorem \ref{lem}, $\Lambda$ is plump in $\Gamma$. 
\end{proof}
 \end{cor} 
 \begin{cor}
\label{3.3}
Let $\Gamma$ be a $C^*$-simple group and let $\Lambda \le \Gamma$ be a non-trivial normal subgroup. If $\Lambda$ contains an element $s$ with amenable centralizer $C_{\Gamma}(s)$ in $\Gamma$, then $\Lambda$ is plump in $\Gamma$. 
\begin{proof} In light of Corollary \ref{3.2}, it is enough to show that $C_{\Gamma}(\Lambda)=\{e\}$. Since $\Lambda$ is normal, so is $C_{\Gamma}(\Lambda)$. Moreover, $C_{\Gamma}(\Lambda) \subset C_{\Gamma}(s)$ for a non-trivial $s \in \Lambda$ is amenable by the assumption. Since $\Gamma$ is $C^*$-simple, it does not contain any non-trivial normal amenable subgroups.  Hence, $C_{\Gamma}(\Lambda)=\{e\}$. 
\end{proof}
\end{cor}
\begin{cor}
\label{finiteindex}
Let $\Gamma$ be a C$^*$-simple group and let $\Lambda$ be a finite index subgroup of $\Gamma$. Then $\Lambda$ is plump in $\Gamma$.
\begin{proof}
Since $\Gamma$ is $C^*$-simple, the action $\Gamma\curvearrowright \partial_F\Gamma$ is free \cite[Theorem 1.1]{BKKO}. It follows from \cite[Chapter-II,Lemma 3.2]{Prox} that this action restricted to $\Lambda$ is strongly proximal. Hence, by Theorem \ref{lem}, $\Lambda$ is plump in $\Gamma$.
\end{proof}
\end{cor}
The above results allow us to use certain free boundary actions to conclude plumpness of subgroups. But, in practice, many natural examples of boundary actions (e.g., $\mathbb{F}_n\curvearrowright \partial\mathbb{F}_n$) are only topologically free. Below, we prove some results which provide us ways to conclude plumpness of subgroups from existence of certain topologically free boundary actions. 

Recall that an action $\Gamma \curvearrowright X$ is topologically free if $X\setminus X^s$ is dense in $X$ for every non-trivial element $s \in \Gamma$, where $X^s=\{x\in X: sx=x\}$.
\begin{lemma}
\label{lemmaplump}
 Let $\Gamma$ be a countable discrete group, let $\Lambda$ be a subgroup of $\Gamma$. Suppose that there exists an action $\Gamma\curvearrowright X$, where $X$ is a compact Hausdorff space, such that for each non-trivial element $s \in \Gamma$, the set $X^s$ of fixed points of $s$, is countable. Further suppose that the action restricted to $\Lambda$ is strongly proximal and that $X$ does not contain any  $\Lambda$-fixed point. Then, $\Lambda$ is plump in $\Gamma$.
\begin{proof} Let $\varphi$ be a state on $C_{\lambda}^*(\Gamma)$. Extend $\varphi$ to a state $\tilde{\varphi}$ on $C(X)\rtimes_r\Gamma$ and let $\tilde{\varphi}|_{C(X)}=d\nu$, where $\nu \in \text{Prob}(\Gamma)$. Since the action restricted to $\Lambda$ is strongly proximal, there are $s_i \in \Lambda$ such that $s_i\nu \to \delta_{x_0}$ in weak$^*$-topology, for some $x_0 \in X$. Now, we claim that $\overline{\Lambda {x_0}}$ is an uncountable set. Let $Y$ be a minimal $\Lambda$-component of $\overline{\Lambda {x_0}}$. If $\overline{\Lambda {x_0}}$ were a countable set, then $Y$ would be a finite set. Since the action restricted to $\Lambda$ is strongly proximal, we must have that $Y$ is a singleton and hence a $\Lambda$-fixed point. This shows that $\overline{\Lambda {x_0}}$ is an uncountable set. Now, since $\cup_{s \ne e, s \in \Gamma}X^s$ is  countable, we can find $y_0 \in \overline{\Lambda {x_0}} \subset \overline{\Lambda\nu}$ with trivial stabilizer. Now, it follows from the proof of \cite[Theorem 4.5]{Haagerup} that $\tau_0 \in \overline{\{s\varphi: s \in \Lambda\}}^{\text{weak}^*}$. Then, by the same theorem, for any finite collection  $t_1,t_2, \ldots, t_n \in \Gamma-\{e\}$ and $\epsilon > 0$, there exist $s_1,s_2,\ldots,s_m \in \Lambda$ such that \[
\left\|\frac{1}{m}\sum_{j=1}^m\lambda_{s_j}\lambda_{t_i}\lambda_{s_j^{-1}}\right\| < \epsilon, \text{ for each $i=1,2,\ldots,n$ }.
\]
\end{proof} 
\end{lemma}

In the case of countable groups $\Gamma$, it turns out that a subgroup $\Lambda$ is plump iff there exists a $C^*$-simple measure on $\Gamma$ which is supported on $\Lambda$.
We recall from \cite{HartKal} that a probability $\mu$ on $\Gamma$ is $C^*$-simple if the canonical trace $\tau_0$ is the unique $\mu$-stationary state on $C_{\lambda}^*(\Gamma)$, that is, $\tau_0$ is the unique state on $C_{\lambda}^*(\Gamma)$ satisfying \[\mu*\tau_0=\sum_{s \in \Gamma}\mu(s)s\tau_0=\tau_0.\]
\begin{prop} Let $\Gamma$ be a countable discrete group. A subgroup $\Lambda \le \Gamma$ is plump in $\Gamma$ iff there exists a $C^*$-simple measure $\mu\in \text{Prob}(\Gamma)$ such that $\text{supp}(\mu) \subseteq \Lambda$. 
\label{prop}
\begin{proof}
Suppose that $\Lambda$ is plump in $\Gamma$. Then, for any finite collection  $t_1,t_2, \ldots, t_n \in \Gamma-\{e\}$ and $\epsilon > 0$, there exist $s_1,s_2,\ldots,s_m \in \Lambda$ such that 
\[
\left\|\frac{1}{m}\sum_{j=1}^m\lambda_{s_j}\lambda_{t_i}\lambda_{s_j^{-1}}\right\| < \epsilon, \text{ for each $i=1,2,\ldots,n$}.
\]
Now arguing similarly, as in the proof of \cite[Theorem 5.1]{HartKal}, one obtains a sequence $\mu_n \in \text{Prob}(\Gamma)$ supported inside $\Lambda$ such that \[\left\|\mu_n*a-\tau_0(a)\mathbf{1}_{C_{\lambda}^*(\Gamma)}\right\| \xrightarrow{n \to \infty}0\quad \forall a \in C_{\lambda}^*(\Gamma).\]
Thus, the $C^*$-simple measure $\mu$ constructed from the sequence $\{\mu_n\}$, as in the proof of \cite[Theorem $5.1$]{HartKal} is also supported inside $\Lambda$. 
 
The converse follows from \cite[Proposition 4.7]{HartKal} and a standard approximation argument.
\end{proof}
\end{prop}
\hfill{}
\subsection{Convergence Groups}
In this section, we use the above results to give some concrete examples where our Theorem \ref{Main Theorem} applies. In particular, we prove Theorem \ref{example}.

An action $\Gamma \curvearrowright X$ is called a convergence action (in this case, $\Gamma$  is called a convergence group), if for every infinite sequence of distinct elements $\gamma _{n}\in \Gamma$, there exist a subsequence $\gamma _{n_{k}}$ and points $a,b\in X$ such that $\gamma _{n_{k}}|_{X\setminus \{a\}}$ converge uniformly on compact subsets to $b$. A non-torsion element $s \in \Gamma$ is called loxodromic if $s$ has exactly two fixed points and is called parabolic if $s$ fixes exactly one point.  A subgroup $\Lambda \le \Gamma$ is elementary if it is finite, or preserves setwise a nonempty subset of $X$ with at most two elements, and called non-elementary otherwise. We refer the reader to \cite{Bow,Fred,Tuk93} for more details on these. 
\begin{prop}
\label{convergenceplump}
Let $\Gamma$ be a non-elementary torsion-free convergence group. Then every non-elementary subgroup $\Lambda \le \Gamma$ is plump in $\Gamma$. 
\begin{proof}
Let $\Gamma \curvearrowright X$ be a convergence action. Since $\Gamma$ is torsion-free, every element in $\Gamma$ is either parabolic or loxodromic (\cite[Theorem 2B]{Tuk93}). Thus, for each non-trivial element $s\in \Gamma$, $|X^s| \le 2$. Since $\Lambda$ is a non-elementary subgroup of $\Gamma$, one sees that the action restricted to $\Lambda$ is strongly proximal (see e.g.,\cite[Example 2]{Oznote}). Since $\Lambda$ is non-elementary, it follows from \cite[Theorem 2S]{Tuk93} that $\overline{\Lambda x}$ is non-trivial for every $x \in X$. The claim now follows from Lemma \ref{lemmaplump}.
\end{proof}
\end{prop}
The above proposition together with Theorem \ref{Main Theorem} imply Theorem \ref{example} for all torsion-free hyperbolic groups. We give an alternative proof below which applies to all hyperbolic groups. 
\hfill{}
\hfill{}

\noindent\\
{\it Proof of Theorem \ref{example}.}\  ~Let $\Gamma$ be a hyperbolic group with trivial amenable radical, let $\Gamma \curvearrowright \mathcal{A}$ be a non-faithful action. Let $\Lambda$ be the kernel of the action $\Gamma \curvearrowright \mathcal{A}$. Since $\Lambda$ is non-amenable, therefore non-elementary, it contains a non-torsion element $s$. We will show that the centraliser $C_{\Gamma}(s)$, of $s$ in $\Gamma$ is amenable, which will imply the theorem by Corollary \ref{3.3}. Let $x_s^+,x_s^- \in \partial\Gamma$ be the points of attraction and repulsion of $s$, i.e. $x_s^+$ and $x_s^{-}$ are fixed by $s$, and $s^nx \xrightarrow{n \to \infty} x_s^{+}$ for all $x \in \partial \Gamma\setminus\{x_s^{-}\}$. We claim that $C_{\Gamma}(s)$ leaves the set $\{x_s^+,x_s^-\}$ invariant. To see this, observe that for any element $t \in C_{\Gamma}(s)$, $tx_s^{\pm}=ts^nx_s^{\pm}=s^ntx_s^{\pm}$. Letting $n \to \infty$, we get that $tx_s^{\pm}=x_s^{\pm}$. Therefore, the kernel of the homomorphism from $C_{\Gamma}(s)$ to $\text{Sym}(\{x_s^+,x_s^-\})$ has index $\le 2$. Let's denote the kernel of this homomorphism by $K$. Since $\Gamma\curvearrowright \partial\Gamma$ is topologically amenable \cite{ADAMS}, the stabilizer $\Gamma_{x_s^+}$ is amenable, hence $K \subset \Gamma_{x_s^+}$ is amenable. Therefore, $C_{\Gamma}(s)$ is amenable. This completes the proof. \qed

The following two examples are general constructions of non-faithful actions for which Theorem \ref{Main Theorem} applies. 

We denote the free group generated by a set $S$ with $\mathbb{F}_{S}$.
 \begin{example} Let $\Gamma$ be a discrete group, let $S$ be a generating set for $\Gamma$ such that $\Gamma \ne \mathbb{F}_S$.  Let $\varphi: \mathbb{F}_{S} \to \Gamma$ be the canonical surjective homomorphism. Let $\alpha: \Gamma \to \rm{Aut}(\mathcal{A})$ be any action of $\Gamma$ on a unital $C^*$-algebra $\mathcal{A}$, let $\tilde{\alpha}=\alpha \circ \varphi$ be the composition action of $\mathbb{F}_S$ on $\mathcal{A}$. Since, $\text{Ker}(\varphi)$ is non-trivial, the action $\mathbb{F}_{S} \curvearrowright \mathcal{A}$ is non-faithful. By Corollary \ref{3.3}, every intermediate $C^*$-subalgebra, $C_{\lambda}^*(\mathbb{F}_{S}) \subseteq \mathcal{B} \subseteq \mathcal{A}\rtimes_{r}\mathbb{F}_{S}$, is of the form $\mathbb{E}(\mathcal{B})\rtimes_{\tilde{\alpha},r}\mathbb{F}_{S}$. 
 \end{example}

\begin{example}
Let $\Gamma$ be a $C^*$-simple group, and let $\Lambda$ be a normal subgroup with trivial centralizer in $\Gamma$. Let $X:=\{0,1\}^{\Gamma/\Lambda}$ and consider the action of $\Gamma$ on $X$ defined by $gf(\bar{s})=f(g^{-1}\bar{s})$, for $f\in X$, $g\in \Gamma$ and $\bar{s} \in \Gamma/\Lambda$. Since $\Lambda$ is plump in $\Gamma$ by Corollary \ref{3.2}, and $\Lambda$ acts trivially on $X$, it follows from Theorem \ref{Main Theorem} that every intermediate $C^*$-subalgebra $\mathcal{B}$, $C_{\lambda}^*(\Gamma) \subseteq \mathcal{B} \subseteq C(X)\rtimes_r \Gamma$, is of the form $\mathbb{E}(\mathcal{B})\rtimes_r\Gamma$. 
\end{example}

\section*{Appendix: Example of a faithful action with a complete description of intermediate algebras}

In this appendix, we give an example of a faithful action of a $C^*$-simple group on a compact Hausorff space $X$, for which every intermediate $C^*$-algebra $\mathcal{B}$, $C_{\lambda}^*(\Gamma)\subseteq \mathcal{B} \subseteq C(X)\rtimes_r\Gamma$ is of the form $\mathcal{A}_1\rtimes_r\Gamma$, where $\mathcal{A}_1$ is a unital $\Gamma$-$C^*$-subalgebra of $C(X)$.
\begin{lemma}
\label{appelemma}
 Let $\Gamma$ be a $C^*$-simple group acting on a compact Hausdorff space $X$ via homeomorphisms. Suppose that the action has the following property:  Given any finite open covering $\{U_i\}$ of $X$, there exists a plump subgroup $\Lambda$ such that $sU_i \subset U_i$ for all $s \in \Lambda$ and every $i$. Then, every intermediate $C^*$-algebra $\mathcal{B}$ with $C_{\lambda}^*(\Gamma) \subseteq \mathcal{B} \subseteq C(X)\rtimes_r\Gamma$ is of the form $\mathcal{A}_1\rtimes_r\Gamma$, where $\mathcal{A}_1$ is a $\Gamma$-$C^*$-subalgebra of $C(X)$. 
\begin{proof} Let $b \in \mathcal{B}$.
Fix $\epsilon>0$. Since $X$ is compact and $\mathbb{E}(b)$ is continuous, we can find finitely many points $\{x_1,\ldots, x_n\}\subset X$ and open neighborhoods containing $x_i$ such that $|\mathbb{E}(b)(x)-\mathbb{E}(b)(x_i)|\leq \epsilon$ whenever $x \in U_{x_i}$ and $X\subseteq \cup_{i}U_{x_i}$. Now let $\Lambda\leq \Gamma$ be the plump subgroup for which the above property holds i.e., $sU_{x_i} \subseteq U_{x_i}$ for all $s \in \Lambda$ and for all $i=1,2,\ldots,n$. Moreover, since $s^{-1}U_{x_i} \subset U_{x_i}$, we get that $sU_{x_i}=U_{x_i}$ forall $i=1,2,\ldots,n$ and for all $s \in \Lambda$. Since $\Lambda$ is plump in $\Gamma$, there exist $s_1,\ldots, s_m\in \Lambda$ such that \[\left\|\frac{1}{m}\sum_{j=1}^m\lambda_{s_j}(b-E(b))\lambda_{s_j^{-1}}\right\|_{\mathbb{B}(\ell^2(\Gamma,\mathcal{H}))}< \epsilon.\] Now, for any $x\in X$, since $X\subseteq \cup_{i=1}^nU_{x_i}$, there exists $1 \le i \le n$ such that $x \in U_{x_i}$. Since $U_{x_i}$ is $\Lambda$-invariant, we see that $s_j^{-1}x \in U_{x_i}$ and hence, by triangle inequality we get that\[\left|\mathbb{E}(b)(s_j^{-1}x)-\mathbb{E}(b)(x)\right|<2\epsilon.\]
Therefore, \[\begin{split}\left\|\frac{1}{m}\sum_{j=1}^m\lambda_{s_j}\mathbb{E}(b)\lambda_{s_j^{-1}}-\mathbb{E}(b)\right\|_{\mathbb{B}(\ell^2(\Gamma,\mathcal{H}))}=&\left\|\frac{1}{m}\sum_{j=1}^m\lambda_{s_j}\mathbb{E}(b)\lambda_{s_j^{-1}}-\mathbb{E}(b)\right\|_{C(X)}\\=&\sup_{x\in X}\left|\frac{1}{m}\sum_{j=1}^m\left(\mathbb{E}(b)(s_j^{-1}x)-\mathbb{E}(b)(x)\right)\right|\\\le&  2\epsilon\\\end{split}\]
By Triangle inequality we get that 
\[\begin{split}\left\|\frac{1}{m}\sum_{j=1}^m\lambda_{s_j}b\lambda_{s_j^{-1}}-\mathbb{E}(b)\right\|_{\mathbb{B}(\ell^2(\Gamma,\mathcal{H}))}\\\le& \left\|\frac{1}{m}\sum_{j=1}^m\lambda_{s_j}(b-\mathbb{E}(b))\lambda_{s_j^{-1}}\right\|_{\mathbb{B}(\ell^2(\Gamma,\mathcal{H}))}\\+&\left\|\frac{1}{m}\sum_{j=1}^m\lambda_{s_j}\mathbb{E}(b)\lambda_{s_j^{-1}}-\mathbb{E}(b)\right\|_{\mathbb{B}(\ell^2(\Gamma,\mathcal{H}))}\\\le& 3\epsilon \\ \end{split}\]
Since $\epsilon>0$ is arbitrary, we get that $\mathbb{E}(b) \in \mathcal{B}$. The assertion now follows from \cite[Proposition 3.4]{Suz17}.
\end{proof}
\end{lemma}
Below we give an example of an action which satisfies the conditions of Lemma \ref{appelemma}. Before we begin the proof, let's briefly recall the notion of odometer actions and some of its properties, which we shall make use of. We refer the reader to \cite{Odo} for more details. 

Let $\Gamma$ be a group and $\{\Gamma_n\}_{n\geq 0}$ be a decreasing sequence of finite-index subgroups. Let $\pi_n: \Gamma/\Gamma_n\to \Gamma/\Gamma_{n-1}$ be the natural quotient map. Consider the inverse limit $X=\lim_n(\Gamma/\Gamma_n, \pi_n)$. We remind that $X$ consists of tuples $(g_0,g_1,g_2,\ldots)\in \prod_{n=0}^{\infty}\Gamma/\Gamma_n$ such that $\pi_n(g_n)=g_{n-1}$ for all $n\geq 1$. The topology on $X$ is generated by the clopen set $\{\{g_n\}\in X: ~g_i=a_i\}$, where $a_i\in \Gamma/\Gamma_i$. 

$\Gamma$ acts on $X$ continuously by left multiplication, i.e. $g\cdot \{h_n\}=\{gh_n\}$, where $g\in \Gamma$ and $\{h_n\}\in X$. This action is called $\Gamma$-odometer or, simply, an odometer when the group $\Gamma$ is clear from the context. If the finite-index subgroups $\{\Gamma_i\}$ are normal, then we call the action an exact $\Gamma$-odometer. Moreover every $\Gamma$-odometer action is minimal and equicontinuous. Also, if the subgroups $\{\Gamma_i\}$ are normal in $\Gamma$, then the above action is free (and hence faithful) if and only if $\cap_{n=0}^{\infty}\Gamma_n=\{e\}$.

\begin{cor} Let $\Gamma$ be a residually finite, $C^*$-simple group. If the action $\Gamma \curvearrowright X$ is a free exact $\Gamma$-odometer, then every intermediate $C^*$-algebra $\mathcal{B}$, $C_{\lambda}^*(\Gamma) \subseteq \mathcal{B}\subseteq C(X)\rtimes_r\Gamma$, is of the form $\mathcal{A}_1\rtimes_r\Gamma$,where $\mathcal{A}_1$ is a $\Gamma$-$C^*$-subalgebra of $\mathcal{A}$.
\begin{proof}
By definition of an exact $\Gamma$-odometer, $X=lim_n(\Gamma/\Gamma_n, \pi_n)$ for some decreasing finite-index normal subgroup $\{\Gamma_i\}$ with $\cap_{n=0}^{\infty}\Gamma_n=\{e\}$. From the definition of the topology on $X$ and since $\cap_{n=0}^{\infty}\Gamma_n=\{e\}$, without loss of generality, after refining the open cover, we may assume the finite open cover to be just the standard one: $\mathcal{C}=\{C_i: 1\leq i\leq [\Gamma: \Gamma_N]\}$ for some large enough $N$, where $C_i=\{\{g_i\}\in X: g_i=a_i\}$ and $\Gamma/\Gamma_N=\{a_i: 1\le i \le [\Gamma:\Gamma_N]\}$. Since $\Gamma_N$ is normal in $\Gamma$, for every $g\in \Gamma_N$, $gC_i=C_i$ for all $i$ as $g$ fixes each coset of $\Gamma/\Gamma_N$. Since $[\Gamma: \Gamma_{N}] < \infty$, $\Gamma_N$ is plump by Corollary \ref{finiteindex}. The assertion now follows from Lemma \ref{appelemma}.
\end{proof}
\end{cor}

\end{document}